\documentclass{article}
\usepackage{graphicx} 
\usepackage{graphicx} 
\usepackage{comment}
\usepackage{amssymb}
\usepackage{amsthm}
\usepackage{extarrows}
\usepackage{calrsfs}
\usepackage{tikz-cd}
\usepackage{mathrsfs}
\usepackage{dynkin-diagrams}
\usepackage{array}
\usepackage[colorlinks = true,citecolor=blue,linkcolor=blue]{hyperref}

\usepackage{todonotes}

\newtheorem{definition}{Definition}[section]
\newtheorem{theorem}[definition]{Theorem}
\newtheorem{lemma}[definition]{Lemma}

\newtheorem{proposition}[definition]{Proposition}

\newtheorem{corollary}[definition]{Corollary}

\newtheorem{remark}[definition]{Remark}

\newcommand{\m}{\mathfrak m}
\newcommand{\kalg}{\mathfrak k}
\newcommand{\Ricop}{\operatorname{Ric}^{\mathrm{op}}}

\usepackage{authblk}
\usepackage{booktabs}
\usepackage{array}
\usepackage{adjustbox}

\newcolumntype{L}{>{\raggedright\arraybackslash}p{0.45\linewidth}}
\newcolumntype{C}{>{\centering\arraybackslash}p{0.45\linewidth}}

\newcommand{\got}{\mathfrak}

\title{
Gray II, Kähler Nullity and Ricci Flow on\\
Generalized Flag Manifolds with Two Isotropy Summands
}

\author{%
  \begin{minipage}[t]{0.31\textwidth}
    \centering
    \textbf{Lino Grama}\par
    \vspace{4pt}
    {\small UNICAMP\par Campinas, Brazil\par}
    {\scriptsize\texttt{lgrama@unicamp.br}\par}
  \end{minipage}%
  \hspace{0.015\textwidth}%
  \begin{minipage}[t]{0.31\textwidth}
    \centering
    \textbf{Ailton R. Oliveira}\par
    \vspace{4pt}
    {\small UEMS\par Dourados, Brazil\par}
    {\scriptsize\texttt{aroliveira@uems.br}\par}
  \end{minipage}%
  \hspace{0.015\textwidth}%
  \begin{minipage}[t]{0.31\textwidth}
    \centering
    \textbf{Francisco V. Oliveira}\par
    \vspace{4pt}
    {\small UFERSA\par Angicos, Brazil\par}
    {\scriptsize\texttt{francisco.vieira@ufersa.edu.br}\par}
  \end{minipage}%
}

\date{}

\begin{document}

\maketitle

\section*{Declarations}

\subsection*{Availability of data and material}
Data sharing is not applicable to this article.

\subsection*{Competing interests}
The authors declare that they have no competing interests.

\subsection*{Funding}
The first author is supported by CNPq grant no.306021/2024-2.

\subsection*{Authors' contributions}
All authors contributed equally to this work.

\subsection*{Acknowledgements}
Not applicable.

\begin{abstract}
We investigate Gray's second curvature identity, Kähler
nullity, and homogeneous Ricci flow on generalized flag
manifolds with two isotropy summands. We prove that every
invariant metric satisfies Gray's second identity for both
invariant almost complex structures. For the non-integrable
structure, the Kähler nullity is trivial except at a
distinguished metric ratio, where it coincides with the
vertical distribution of a homogeneous fibration. This
distribution is integrable and has totally geodesic leaves.
For the integrable structure, the Kähler nullity exhibits
a different behavior and need not be integrable away from
the Kähler locus. Finally, we reduce the normalized
homogeneous Ricci flow to an autonomous equation for the
ratio of the metric parameters. The resulting dynamics
identify two invariant Einstein directions and show that
the Kähler--Einstein direction is attracting, while the
other is repelling.
\end{abstract}

\noindent\textbf{Keywords:}
Generalized flag manifolds; Gray curvature identities;
Kähler nullity; invariant almost complex structures;
homogeneous Ricci flow; Einstein metrics.

\section{Introduction}

Almost Hermitian geometry provides a natural framework for studying the interaction between Riemannian metrics and almost complex structures.  In this setting, curvature identities play a central role and have been studied by several authors; see, for instance, \cite{alex,davi,navi}. Gray~\cite{gray} observed that the curvature tensor of a Kähler manifold satisfies the identity
$R_{WXJYJZ}=R_{WXYZ}$.
Motivated by this property, he investigated which curvature identities hold for more general classes of almost Hermitian manifolds. In the same work, he introduced the following identities:
\begin{align*}
(G_{1}):& \quad R_{WXYZ}=R_{WXJYJZ},\\
(G_{2}):& \quad R_{WXYZ}-R_{JWJXYZ}
    =R_{JWXJYZ}+R_{JWXYJZ}, \label{G2}\\
(G_{3}):& \quad R_{WXYZ}=R_{JWJXJYJZ}.
\end{align*}

A fundamental result of Gray \cite{gray} states that
the K\"ahler nullity distribution of a
quasi-K\"ahler manifold satisfying the
second Gray curvature identity is
integrable wherever its rank is constant,
and its integral manifolds are K\"ahler
submanifolds. We recall this result
in Theorem~2.12.

It is well known that a quasi-Kähler manifold endowed with an integrable almost complex structure is Kähler \cite{GrayHervella1980}. We therefore focus primarily on the
non-integrable case. The integrable case will be discussed in the section \ref{integrable}.

Homogeneous spaces, and in particular generalized flag manifolds, provide a convenient setting in which these questions can be investigated explicitly. Their invariant geometric structures can be described in terms of Lie-theoretic data, while invariant Riemannian metrics and almost complex structures are determined by a finite number of parameters. This makes generalized flag manifolds a natural class of examples for studying curvature identities beyond the Kähler setting.

Throughout the statements of our main results, we consider
the class of generalized flag manifolds with two isotropy
summands described in Section~3.
Our main observation is that Gray's second curvature identity imposes no restriction on the metric parameters: 
\begin{theorem}\label{prop:all-gray}
Let $G/K$ be a generalized flag manifold with two isotropy summands and let
$J=(+,-)$ be the invariant almost complex structure. Every invariant Hermitian metric \eqref{eq:metric-two} satisfies
Gray's second curvature identity.  Thus, within the invariant metric family,
\begin{equation}\label{eq:gray-whole-family}
\{g_{\lambda}:G_2(g_{\lambda},J)=0\}
 =\{g_{\lambda}:\lambda_1>0,\ \lambda_2>0\}.
\end{equation}
\end{theorem}

Theorem~1.1 is proved in Theorem~3.2 using the root-space decomposition and the possible degree configurations of complementary roots. As a consequence, Gray's third identity is also satisfied throughout the invariant metric family.
 We remark that Theorem~1.1 was previously obtained in \cite{AAD} using methods different from those employed here.

It is important to distinguish K\"ahler nullity from the
classical nullity of the Riemannian curvature tensor.
The latter has been extensively studied in homogeneous
geometry, notably by Di Scala, Olmos and Vittone
\cite{DiScalaOlmosVittone2022}.
Unlike the classical nullity of the Riemannian
curvature tensor, Kähler nullity is defined
through the covariant derivative of the
almost complex structure and may determine
a nontrivial homogeneous distribution.

In our context, assuming that
\[
[\mathfrak m_1,\mathfrak m_1]_{\mathfrak m_2}\neq 0,
\]
we show that the Kähler nullity is trivial for $(\lambda_2\neq 2\lambda_1)$, whereas for $(\lambda_2=2\lambda_1)$ it coincides with the invariant distribution determined by $(\mathfrak m_2)$. Independently of the almost complex structure, this distribution is integrable and its leaves are totally geodesic. These facts are summarized in the following two results.

\begin{theorem}
Let
$
M=G/K
$
be a compact homogeneous flag manifold with two isotropy summands
$
\got{m}=\got{m}_1\oplus\got{m}_2
$
and invariant metric
$
g.
$
Then the distribution
\[
\mathcal{D}_2=G\times_K\got{m}_2
\]
is integrable and its leaves are totally geodesic.
\end{theorem}

\begin{theorem}
Let $M=G/K$ be a flag manifold with two isotropy summands and $J=(+,-)$ be the non-integrable almost complex structure. If $[\mathfrak{m}_{1},\mathfrak{m}_{1}]_{\mathfrak{m}_{2}}\ne0$, then,
$
\mathcal{K}=\left\{\begin{array}{cc}
   0,  & \lambda_{2} \neq 2\lambda_{1} \\
   \mathcal{D}_{2},  & \lambda_{2}= 2\lambda_{1}
\end{array}\right..   
$
\end{theorem}

In the example of $(G_2/U(2))$, it is geometrically identified with the vertical distribution of the twistor fibration
\[
G_2/U(2)\longrightarrow G_2/SO(4),
\]
whose leaves are copies of $(\mathbb{CP}^1)$.

We also examine the integrable invariant complex structure ($J=(+,+)$). In this case, Gray's second identity is again satisfied by the entire family of invariant metrics, but the behavior of the Kähler nullity is substantially different. Away from the Kähler metric, the nullity is determined by the first isotropy summand, whereas at the distinguished ratio ($\lambda_2=2\lambda_1$) the structure becomes Kähler and the nullity coincides with the whole tangent bundle. This comparison emphasizes that Gray's second identity alone does not force the Kähler nullity distribution to be integrable and highlights the role of the quasi-Kähler hypothesis in Gray's integrability Theorem.

Beyond their intrinsic geometric interest, homogeneous
Ricci flows also arise naturally in the study of
renormalization group equations for nonlinear sigma
models. At one-loop order, the renormalization group
evolution of the target-space metric is governed by
the Ricci tensor, up to conventional rescaling and
reparameterization \cite{Friedan1980}.
In this context, the restriction to invariant metrics
provides a finite-dimensional setting in which the
evolution of relative metric parameters can be
analyzed explicitly. Although our results are purely
geometric, the dynamical behavior of the invariant
Einstein directions may be viewed in connection with
this broader geometric framework.

Finally, we investigate the behavior of these structures under the homogeneous Ricci flow. Recall from \cite{ham-ricci}that the Ricci ﬂow equation is defined by
\begin{equation}\label{def-flow}
\dfrac{\partial g(t)}{\partial t}=-2Ric(g(t)),
\end{equation}
where Ric(g) is the Ricci tensor of the Riemannian metric g. The solution of this equation, The resulting one-parameter family of metrics on
$M$ is called the Ricci flow. Intuitively, this is the heat
equation for the metric g. The Ricci ﬂow equation \ref{def-flow} for an arbitrary manifold is a nonlinear system of PDEs. However,
when restricted to the set of invariant metrics, system \ref{def-flow} reduces to an autonomous nonlinear
system of ODEs. Because of this, it is natural to proceed with a qualitative study, using tools from the theory of dynamical systems.

Since the flow preserves the family of invariant metrics, and Gray's second identity holds on the whole positive parameter space, the Gray condition is preserved along every invariant Ricci flow trajectory. Introducing the ratio
\[
u=\frac{\lambda_2}{\lambda_1},
\]
we reduce the normalized homogeneous Ricci flow to a one-dimensional autonomous equation. The two equilibrium values correspond to the two invariant Einstein directions: one is the Kähler-Einstein ratio ($u=2$), while the other corresponds to the second invariant Einstein metric. The resulting phase portrait shows that the Kähler-Einstein direction is attracting in the projectivized space of invariant metrics, whereas the second Einstein direction is repelling.

The paper is organized as follows. In Section 2, we recall the basic facts concerning invariant metrics and almost complex structures on generalized flag manifolds, quasi-Kähler structures, Gray's curvature identities, and Kähler nullity. In Section 3, we specialize to flag manifolds with two isotropy summands and study Gray's second identity, invariant distributions, and the Kähler nullity, with particular attention to ($G_2/U(2)$) and $(\mathbb{CP}^3=Sp(2)/(Sp(1)\times U(1)))$. In Section 4, we analyze the homogeneous Ricci flow and derive the evolution equation for the ratio of the metric parameters, describing the corresponding Einstein directions and phase portrait.

\section{Invariant geometry of flag manifolds}\label{sec:prelim}

In this section we recall some well known results about the almost Hermitian geometry of generalized flag manifolds, from a Lie theoretical point of view. See for example \cite{gau1}, \cite{kobaiashi}, \cite{livrosm},\cite{besse}, \cite{sr}, \cite{neiton-luciana} and \cite{sn}.  

\subsection{Generalized flag manifolds} \label{sec-flag-review}
Let $\got{g}^\mathbb{C}$ be a complex semi-simple Lie algebra. Fix a Cartan subalgebra $\got{h}$ of $\got{g}^\mathbb{C}$ and denote by $\Pi$ the corresponding root system with respect to the pair $(\got{g}^\mathbb{C},\got{h})$. We have the following root space decomposition
$$
\got{g}^\mathbb{C}=\got{h}\oplus\sum_{\alpha\in \Pi}\got{g}_{\alpha},
$$ 
where $\got{g}_{\alpha}=\{X\in\got{g}^{\mathbb{C}}:\forall H\in \got{h},[H,X]=\alpha(H)X \}$ denote the corresponding $1$-dimensional (complex) root space. 

The Cartan-Killing form is defined by  
$$
\langle X,Y\rangle=tr(ad(X)ad(Y))
$$ 
and its restriction to $\got{h}$ is non-degenerate. Given a root $\alpha\in \got{h}^{*}$ we define $H_{\alpha}$ by $\alpha(\cdot)=\langle H_{\alpha},\cdot\rangle$. Moreover we denote $\got{h}_{\mathbb{R}}=span_{\mathbb{R}}\{H_{\alpha}:\alpha\in \Pi\}$ and $\got{h}_{\mathbb{R}}^{*}$ being the real subspace of $\got{h}^{*}$ spanned by the roots.

We will fix a Weyl basis of $\got{g}^\mathbb{C}$ as follows:
$$
X_{\alpha}\in\got{g}_{\alpha}, \ \mbox{ such that } \ \langle X_{\alpha},X_{-\alpha}\rangle=1 \ \mbox{and} \ [X_{\alpha},X_{\beta}]=m_{\alpha,\beta}X_{\alpha+\beta},
$$ 
with $m_{\alpha,\beta}\in\mathbb{R}$, $m_{-\alpha,-\beta}=-m_{\alpha,\beta}$, $[X_{\alpha},X_{-\alpha}]=H_{\alpha}$ and $m_{\alpha,\beta}=0$ if $\alpha+\beta$ is not a root. 

Let $\Pi^{+}\subset \Pi$ be a choice of positive roots, $\Sigma$ be the corresponding system of simple roots and $\Theta$ be a subset of $\Sigma$. We will use the following notation: $\langle\Theta\rangle$ is the set of roots spanned by $\Theta$, $\Pi_{M}=\Pi\setminus\langle{\Theta}\rangle$ be the set of complementary roots and $\Pi_{M}^{+}$ be the set of complementary positive roots.

Let 
$$
\got{p}_{\Theta}=\got{h}\oplus
\sum_{\alpha\in\langle\Theta\rangle^{+}}\got{g}_{\alpha}\oplus
\sum_{\alpha\in\langle\Theta\rangle^{+}}\got{g}_{-\alpha}\oplus
\sum_{\beta\in \Pi_{M}^{+}}\got{g}_{\beta}
$$
be a parabolic sub-algebra of $\got{g}^\mathbb{C}$ determined by $\Theta$.

The generalized flag manifold $\mathbb{F}_{\Theta}$ is the homogeneous space 
$$
\mathbb{F}_{\Theta}=G^\mathbb{C}/P_{\Theta},
$$
where $G^\mathbb{C}$ is a complex connected Lie group with Lie algebra $\got{g}^\mathbb{C}$ and $P_{\Theta}$ is the normalizer of $\got{p}_{\Theta}$ in $G^\mathbb{C}$.

Let $\got{g}$ be the compact real form of $\got{g}^\mathbb{C}$. 
We have
$$
\got{g}=span_{\mathbb{R}}\{i\got{h}_{\mathbb{R}},A_{\alpha},iS_{\alpha}; \alpha\in \Pi\},
$$ 
where $A_{\alpha}=X_{\alpha}-X_{-\alpha}$ and $S_{\alpha}=X_{\alpha}+X_{-\alpha}$. We remark that the Lie algebra $\got{g}$ is semi-simple.

Denote by $G$ the compact real form of $G^\mathbb{C}$ with $Lie(G)=\got{g}$ and let $\got{k}_{\Theta}$ the Lie algebra of  $K_{\Theta}:=P_{\Theta}\cap G$. It is well known that the Lie group $K_{\Theta}\subset G$ is a centralizer of a torus. 

We have 
$$
\got{k}_{\Theta}^{\mathbb{C}}=\got{h}\oplus
\sum_{\alpha\in\langle\Theta\rangle^{+}}\got{g}_{\alpha}\oplus
\sum_{\alpha\in\langle\Theta\rangle^{+}}\got{g}_{-\alpha},
$$
where $\got{k}_{\Theta}^{\mathbb{C}}$ denotes the complexification of the real Lie algebra $\got{k}_{\Theta}=\got{g}\cap \got{p}_{\Theta}$. 

The Lie group $G$ also acts transitively on $\mathbb{F}_{\Theta}$ and we have 
$$
\mathbb{F}_{\Theta}=G^\mathbb{C}/P_{\Theta}=G/(P_{\Theta}\cap G)=G/K_{\Theta}.
$$

When $\Theta=\emptyset$ the parabolic sub-algebra is a Borel sub-algebra of $\got{g}^\mathbb{C}$ and $\mathbb{F}=G^\mathbb{C}/P=G/T$ is called {\it full flag manifold}, where $T=P\cap G$ is a maximal torus of $G$. 
\bigskip

Recall that the flag manifold $\mathbb{F}_{\Theta}=G/K_{\Theta}$ is a reductive homogeneous space, that is, there exists a subspace $\got{m}$ of $\got{g}$ such that 
$$
\got{g}=\got{k}_{\Theta}\oplus \got{m} \ \ \mbox{and} \ \ 
Ad(k)\got{m}\subseteq \got{m}, \forall k\in K_{\Theta},
$$
and one can identify the tangent space $T_{x_{0}}\mathbb{F}_{\Theta}$ with $\got{m}$, where $x_{0}=eK_{\Theta}$ is the origin of the $\mathbb{F}_{\Theta}$ (trivial coset).

Let us give a description of the tangent space $\got{m}$ in terms of the Lie algebra structure of $\got{g}$ as follows:
$$
\got{g}=\got{k}_{\Theta}\oplus\sum_{\beta\in \Pi_{M}^{+}}{\got{u}_{\beta}},
$$
with $\got{u}_{\beta}=\got{g}\cap(\got{g}_{\beta}\oplus \got{g}_{-\beta})$ and for each root $\beta \in \Pi_{M}$, $\got{u}_{\beta}$ has real dimension two and it is spanned by $A_{\beta}$ and $\sqrt{-1}S_{\beta}$, and we have the following identification $\got{m}=\sum_{\beta\in \Pi^{+}_{M}}\got{u}_{\beta}$.

\

An important tool to study invariant tensors on homogeneous space is the isotropy representation. We will restrict ourselves to the situation of flag manifolds. In this case, since the flag manifolds are reductive homogeneous spaces it is well know that the isotropy representation is equivalent to the following representation
$$Ad(k)\mid_{\got{m}}:\got{m}\longrightarrow\got{m}.$$

The isotropy representation decomposes $\got{m}$ into irreducible components, that is, 
$$
\got{m}=\got{m}_{1}\oplus\got{m}_{2}\oplus\cdots\oplus\got{m}_{n},
$$
where each component $\got{m}_{i}$ satisfies $Ad(K_{\Theta})(\got{m}_{i})\subset \got{m}_{i}$. 
Moreover, each component $\got{m}_{i}$ is irreducible, that is, the only invariant sub-spaces of $\got{m}_{i}$ by $Ad(K_{\Theta})|_{\got{m}_{i}}$ are the trivial sub-spaces. We will call the sub-spaces $\got{m}_{i}$ by {\it isotropy summands} of the isotropy representation. 

\begin{remark}
From now on we will omit the symbol $\Theta$ whenever there is no risk of confusion. We will denote a flag manifold just by $\mathbb{F}=G/K$.
\end{remark}

\subsection{Invariant Metrics} 
Let us denote by $\langle\cdot,\cdot\rangle$ the Cartan-Killing form of $\got{g}$. For each $Ad(K)-$invariant inner product $B$ on $\mathfrak{m}$, there exists a unique $\langle\cdot,\cdot\rangle-$self-adjoint, positive operator $\Lambda :\mathfrak{m}\longrightarrow\mathfrak{m}$ commuting with $Ad(k)\left|_{\mathfrak{m}}\right.$ for all $k\in K$ such that 
$$B(X,Y)=Q(\Lambda X,Y), \,\, X,Y\in\mathfrak{m},$$
where $Q$ is the negative of the Cartan-Killing of $\mathfrak{g}$.

Therefore an invariant Riemannian metric $g$ on $\mathbb{F}$ is completely determined by the invariant inner product $B$, and the inner product is determined by $\Lambda$. 

The vectors $A_{\alpha}$, $\sqrt{-1}S_{\alpha}$, $\alpha\in\Pi_{M}$,  are eigenvectors of $\Lambda$ associated to the same eigenvalue $\lambda_\alpha$. 

The invariant inner product $B$ admits a natural extension to a symmetric bi-linear form on  $\got{m}^{\mathbb{C}}$. On the complexified tangent space we have $\Lambda(X_{\alpha})=\lambda_{\alpha}X_{\alpha}$ with  $\lambda_{\alpha}>0$ and $\lambda_{-\alpha}=\lambda_{\alpha}$. 

\

{\bf Notation:} Throughout the paper, we identify an invariant
metric with its defining operator $\Lambda$. We also denote the invariant metric $g$ just by a $n$-tuple of positive numbers $(\lambda_1,\ldots, \lambda_n)$ representing the eigenvalues of the operator $\Lambda$ and is parameterized by the number of irreducible components.   \qed

\subsection{Invariant Almost Complex Structures}

\begin{definition}
An almost complex structure on a smooth
manifold $M$ is a smooth section
\[
J\in\Gamma(\operatorname{End}(TM))
\]
satisfying
\[
J^2=-\operatorname{Id}_{TM}.
\]
\end{definition}

\begin{definition}
A $G$-invariant almost complex structure $J$ on $\mathbb{F}=G/K$ is an almost complex structure that satisfies  
$$
J_{ux}=dE_{u}J_{x}dE_{u^{-1}},\ \mbox{for all } \ u\in G,
$$ 
where $dE_{u}:T(G/K) \to T(G/K)$ denotes the differential of the left translation by $u$, that is, for all $X\in T_{x}(G/K)$ we have
$$
dE_{u}J_{x}X=J_{ux}dE_{u}X.
$$
\end{definition}

The following result allows us to describe invariant almost complex structures on flag manifolds in terms of complex structure in a simple vector space, namely the tangent space at the origin (trivial coset) of the homogeneous space. 

\begin{proposition} \label{prop-corresp-iacs}
There is a one-to-one correspondence between a $G$-invariant almost complex structure $J$ and a linear endomorphism $J_{x_{0}}: T_{x_{0}}\mathbb{F} \to T_{x_{0}}\mathbb{F}$ satisfying $J_{x_{0}}^2=-Id$ and commute with the isotropy representation, that is, 
$$
Ad^{G/K}(k)J_{x_{0}}=J_{x_{0}}Ad^{G/K}(k), \ \mbox{for all } \ k\in K. 
$$
\end{proposition}

An interesting consequence of the Proposition \ref{prop-corresp-iacs} is that $J(\got{g}_{\alpha})=\got{g}_{\alpha}$, where $\got{g}_\alpha$ is the root space associated to the root $\alpha \in \Pi$. The eigenvalues of $J$ are $\pm \sqrt{-1}$ and the eigenvectors on $\got{m}^{\mathbb{C}}$ are  $X_{\alpha}$, $\alpha\in\Pi$. Therefore $J(X_{\alpha})= \varepsilon_{\alpha} \sqrt{-1} X_{\alpha}$, with $\varepsilon_{\alpha}=\pm 1$ and  $\varepsilon_{\alpha}=-\varepsilon_{-\alpha}$. 

Summarizing, an invariant almost complex structure is determined by a choice of signs on the irreducible components of the complexified isotropy representation (equivalently, on the corresponding \(t\)-roots), subject to \(\varepsilon_{-\alpha}=-\varepsilon_\alpha\).


\begin{definition}[Definition 2.2 of \cite{sn}]
Let $J$ be an invariant almost complex structure on the flag manifold $G/K$. The triple of roots $\alpha, \beta, \gamma\in \Pi_M$, with  $\alpha+\beta+\gamma=0$ is said to be a $(0,3)$-triple if $\varepsilon_\alpha=\varepsilon_\beta=\varepsilon_\gamma$. It is a $(1,2)$-triple otherwise.
\end{definition}

\subsection{Invariant Almost Hermitian Structures}

\begin{definition}
     An Invariant Almost Hermitian Structure  is composed of a pair $(g,J)$ with
$J$ an invariant almost complex structure and $g$ an invariant Riemannian metric.
\end{definition}

\begin{definition}
    The Kähler 2-Form or fundamental 2-form, $\omega$, is defined by $\omega(X,Y)=g(JX,Y)$.
\end{definition}

In this work we are interested in $(1,2)$-symplectic (or Quasi-Kähler) flag manifolds. 
Recall that an almost Hermitian manifold is said to be (1,2)-symplectic (or quasi-Kähler) if
$$
d\omega(X,Y,Z)=0
$$
when one of the vectors X; Y; Z is of type (1,0) and the other two are
of type (0,1).

Another characterization of such structures is given in the next proposition.

\begin{proposition}[Proposition 2.3 of \cite{sn}]
The invariant pair ($J=\{\varepsilon_{\alpha}\}; \Lambda=\{\lambda_{\alpha}\})$ is $(1,2)$-symplectic if and
only if
$$
\varepsilon_{\alpha}\lambda_{\alpha}+\varepsilon_{\beta}\lambda_{\beta}+\varepsilon_{\gamma}\lambda_{\gamma}=0
$$
for every $(1,2)$-triple.
\end{proposition}

A characterization using the Riemannian connection is given below.
\begin{lemma}[\cite{gau1}]
The pair is $(1,2)$-symplectic if and only if $D_{J\cdot}J=-J\circ D\cdot J$.
\end{lemma}

\subsection{Second Gray Identity}

Consider the following identities, given in \cite{gray}, for the curvature operator of an almost Hermitian manifold $M$: for $W,X,Y,Z \in \mathcal{H}(M)$,
\begin{align}   
(G_{1}):& \ \ R_{WXYZ}=R_{WXJYJZ}\\
(G_{2}):& \ \ R_{WXYZ}-R_{JWJXYZ}=R_{JWXJYZ}+R_{JWXYJZ} \label{G2}\\
(G_{3}):& \ \ R_{WXYZ}=R_{JWJXJYJZ}.
\end{align}
Here $R_{XYZW}$ means $R(X,Y,Z,W)$.
Sometimes we will use
\begin{eqnarray*}
G_{2}(X,Y,Z,W)&=&R(X,Y,Z,W)
 -R(JX,JY,Z,W)\\
&-&R(JX,Y,JZ,W)
 -R(JX,Y,Z,JW).    
\end{eqnarray*}
In this context we want to understand when $G_{2}=0$.

For a given class $\mathcal{L}$ of almost Hermitian manifolds, let $\mathcal{L}_{i}$ be the
subclass of manifolds whose curvature operator satisfies identity (i). It follows that 
$$\mathcal{L}_{1}\subseteq \mathcal{L}_{2} \subseteq \mathcal{L}_{3} \subseteq \mathcal{L}.$$

Kähler manifolds satisfy all identities $G_{1}$, $G_{2}$ and $G_{3}$ and Nearly Kähler manifolds satisfy both $G_{2}$ and $G_{3}$. See \cite{gray}, \cite{gray2} and \cite{sr}.

\begin{definition}[\cite{gray}]
Let $M$ be an almost Hermitian manifold. For each $p\in M$ we set
\[
\mathcal K_p
=
\{X\in T_pM:
(\nabla_XJ)Y=0
\text{ for every }Y\in T_pM
\}.
\]
The family \(p\mapsto\mathcal K_p\) is called the Kähler nullity. On an open set where its dimension is constant, it defines a distribution.
\end{definition}

If $M$ is a Kähler manifold, then $\mathcal{K}(p)=T_{p}M$, for all $p$, since $DJ=0$. What can be said about the Kähler nullity of
non-Kähler almost Hermitian manifolds?

The Nullity Kähler can be calculated using Lemma below. 

\begin{lemma}[\cite{sigma}]\label{nulidade} 
The covariant derivative of the K\"ahler form $\omega$ with respect to the Levi-Civita connection $D$ on a flag manifold $G/K$ is given by
$$
g((D_{X_{\alpha}}J)X_{\beta},X_{\gamma})
\left\{
\begin{array}{ll}
-\dfrac{\sqrt{-1}}{2} m_{\alpha,\beta}(\varepsilon_{\beta}+\varepsilon_{\gamma})(-\lambda_{\alpha}+\lambda_{\beta}+\lambda_{\gamma}) &, \mbox{ if }  \alpha+\beta+\gamma = 0, \\
0 &, \mbox{ otherwise }.
\end{array}
\right.
$$
\end{lemma}

We will use the notation $D\omega(X_{\alpha},X_{\beta},X_{\gamma}):=g((D_{X_{\alpha}}J)X_{\beta},X_{\gamma})$.

Now we will analyze the following result and its geometric applications.

\begin{theorem}[Theorem 7.1 of \cite{gray}]
Let $M\in \mathcal{Q}\mathcal{K}$ (Quasi-Kähler manifold) satisfy the curvature identity ($G_{2}$), namely $R_{WXYZ}-R_{WXJYJZ}=R_{WJXYJZ}+R_{JWXYJZ}$ for $W,X,Y,Z\in \mathcal{H}(M)$. Then on any open subset of $M$ on which dim $\mathcal{K}(m)$ is constant, the distribution $m \rightarrow \mathcal{K}(m)$ is integrable. Furthermore any integral manifold of $m \rightarrow \mathcal{K}(m)$ is a Kähler submanifold of $M$.    
\end{theorem}

It was proved that $m \rightarrow \mathcal{K}(m)$ is integrable provided $M$ is a Nearly Kähler manifold. The same result holds when $M$ is a quasi-Kähler manifold satisfying curvature identity $G_{1}$.

\section{Flag manifolds with two isotropy summands }

The classification of generalized flag manifolds with two isotropy summnands is provided in \cite{AC2011}.
Throughout this section, we assume that $G$ is a compact,
connected, simple Lie group and that $K$ is the centralizer
of a torus determined by a single painted simple root
$\alpha_p$ whose coefficient in the highest root is two,
that is, $c_p=2$.

Under these assumptions, the isotropy representation
decomposes into exactly two irreducible, pairwise
inequivalent summands,
\[
    \mathfrak{m}=\mathfrak{m}_1\oplus\mathfrak{m}_2.
\]

In this section, we investigate Gray's second curvature
identity, K\"ahler nullity, and invariant distributions
on the class of generalized flag manifolds $G/K$
described above. In particular, we study the interplay
between invariant almost Hermitian structures,
the isotropy decomposition
$\mathfrak{m}=\mathfrak{m}_1\oplus\mathfrak{m}_2$,
and the geometric properties of the associated
homogeneous distributions. Our analysis builds on
the foundational results of Gray \cite{gray}
and the homogeneous geometric framework discussed
in \cite{gorb,kobaiashi}.

The classification of this family of homogeneous space is given as follows: recall the {\em height} (or degree) of a simple root $\alpha_p$ is the coefficient of $\alpha_{p}$ in the highest root, that is, if $$
\mu=c_1{\alpha_1}+\ldots + c_p\alpha_p+\ldots +c_n\alpha_n
$$ is the highest root then $ht (\alpha_p)=c_p$. 

A generalized flag manifold has two isotropy components if, and only if, it is parameterized by the set $\Theta\subset\Sigma$ such that $\Theta=\Sigma \setminus\{ \alpha_p\}$ with $ht(\alpha_p)=2$. See \cite{AC2011} for further details. The classification is understood up to
local isomorphism of homogeneous spaces; in
particular, finite coverings and central
quotients do not affect the corresponding
isotropy decompositions or the invariant
geometric calculations considered below.


\begin{table}[!htb] 

\centering

\caption{Flag manifolds with two isotropic summands, \cite{AC2011}} 

\begin{tabular}{|c|c|} 

\hline 

$G$ simple & Flag manifold $G/K$ with $\got{m}=\got{m}_{1}\oplus \got{m}_{2}$\\ 

\hline
\hline
$B_{l}$ & $SO(2l+1)/U(p)\times SO(2(l-p)+1) \ (2\leq p \leq l)$\\
\hline
$C_{l}$ & $Sp(l)/U(p)\times Sp(l-p) \ (1\leq p \leq l-1)$\\
\hline
$D_{l}$ & $SO(2l)/U(p)\times SO(2(l-p)) \ (2\leq p \leq l-2)$\\
\hline
$G_{2}$ & $G_{2}/U(2)$ ($U(2)$ corresponding to the short simple
root of $G_{2}$)\\
\hline
$F_{4}$ & $F_{4}/SO(7)\times U(1)$\\
        & $F_{4}/Sp(3)\times U(1)$\\
\hline
$E_{6}$ & $E_{6}/SU(6)\times U(1)$\\
        & $E_{6}/SU(2)\times SU(5)\times U(1)$\\
\hline
$E_{7}$ & $E_{7}/SU(7)\times U(1)$\\
        & $E_{7}/SU(2)\times SO(10)\times U(1)$\\
        & $E_{7}/SO(12)\times U(1)$\\
\hline
$E_{8}$ & $E_{8}/E_{7}\times U(1)$\\
        & $E_{8}/SO(14)\times U(1)$\\
\hline
\end{tabular}
\label{tab1}
\end{table}

According to \cite{AC2011}, the isotropy summands $\got{m}_{1}$ and $\got{m}_{2}$ satisfy the following inclusions: 
$$
[\mathfrak{m}_{1},\mathfrak{m}_{1}]\subset \mathfrak{k}\oplus\mathfrak{m}_{2}, \ \ \ \ [\mathfrak{m}_{1},\mathfrak{m}_{2}]\subset \mathfrak{m}_{1} \ \ \ \ \mbox{and} \ \ \ \  [\mathfrak{m}_{2},\mathfrak{m}_{2}]\subset \mathfrak{k}.
$$


In this case, an invariant metric on $G/K$ is given by
\begin{equation}\label{eq:metric-two}
\m=\m_1\oplus\m_2,\qquad
 g_{\lambda}=\lambda_1Q|_{\m_1}+\lambda_2Q|_{\m_2},\qquad
 \lambda_1,\lambda_2>0.
\end{equation}

Moreover, $G/K$ admits $2$ invariant almost complex structures, up to conjugation:
$$
\begin{array}{ccccccc}
J_{1}=(+,+)&&&\mbox{and}&&&J_{2}=(+,-).
\end{array}
$$



\begin{lemma}[\cite{sigma}]\label{QK}
Let $G/K$ be a generalized flag manifold with two isotropy summands. If
$J=(+,-)$ is the invariant almost complex structure and $g$ is an invariant metric, then the pair $(g,J)$ is quasi-Kähler.
\end{lemma}

The following Theorem characterizes all metrics that satisfy the second Gray identity. Assume that  $G_2(g_{\lambda},J)=0$ 
means the metric $g_{\lambda}$ and almost complex structure $J$ that satisfy the $G_{2}$ identity.

\begin{theorem}\label{prop:all-gray}
Let $G/K$ be a generalized flag manifold with two isotropy summands and let
$J=(+,-)$ be the invariant almost complex structure. Every invariant Hermitian metric \eqref{eq:metric-two} satisfies
Gray's second curvature identity.  Thus, within the invariant metric family,
\begin{equation}\label{eq:gray-whole-family}
\{g_{\lambda}:G_2(g_{\lambda},J)=0\}
 =\{g_{\lambda}:\lambda_1>0,\ \lambda_2>0\}.
\end{equation}
\end{theorem}

\begin{proof}
The manifolds $G/K$ have complementary root degrees (deg) $1$ and $2$.
Consequently the degrees of complexified tangent root vectors belong to
$\{1,-1,2,-2\}$.  Write $JX_\alpha=i\varepsilon_\alpha X_\alpha$.  For
$J=(+,-)$ their signs are
\begin{equation}\label{eq:degree-signs}
\begin{array}{c|rrrr}
\deg\alpha&1&-1&2&-2\\ \hline
\varepsilon_\alpha&+1&-1&-1&+1.
\end{array}
\end{equation}
Since the curvature tensor is invariant under the isotropy torus,
$R(X_\alpha,X_\beta,X_\gamma,X_\delta)$ vanishes unless
$\alpha+\beta+\gamma+\delta=0$.  In particular, the four degrees must sum
to zero.  Up to permutation, the possibilities are exactly
\begin{equation}\label{eq:four-degrees}
(1,1,-1,-1),\qquad (2,2,-2,-2),\qquad (1,-1,2,-2).
\end{equation}
These are all degree multisets with zero sum; for each of them, every
ordering has two positive and two negative \(J\)-signs.

Each row of \eqref{eq:four-degrees} contains two $+1$ and two $-1$ signs
from \eqref{eq:degree-signs}.  With the convention for the second Gray identity from (\ref{G2}), its defect on a quadruple of root vectors is
\begin{align*}
G_2(X_\alpha,X_\beta,X_\gamma,X_\delta)
&=R(X_\alpha,X_\beta,X_\gamma,X_\delta)
 -R(JX_\alpha,JX_\beta,X_\gamma,X_\delta)\\
&\quad-R(JX_\alpha,X_\beta,JX_\gamma,X_\delta)
 -R(JX_\alpha,X_\beta,X_\gamma,JX_\delta)\\
&=\bigl(1+\varepsilon_\alpha\varepsilon_\beta
 +\varepsilon_\alpha\varepsilon_\gamma
 +\varepsilon_\alpha\varepsilon_\delta\bigr)
 R(X_\alpha,X_\beta,X_\gamma,X_\delta).
\end{align*}
The coefficient is zero whenever there are two signs of each kind.  For every
other quadruple the curvature component itself vanishes by torus invariance.
Therefore $G_2=0$ for every $\lambda_1,\lambda_2>0$.
\end{proof}


\begin{remark}\label{rem:gray-not-ratio}
The relation $\lambda_2=2\lambda_1$ does \emph{not} characterize Gray II
for $J=(+,-)$.  It singles out the invariant K\"ahler--Einstein metric for
the integrable invariant complex  $J=(+,+)$, which is a different geometric
condition.  In particular, the proof of Theorem~\ref{prop:all-gray}
uses neither the nearly K\"ahler condition nor a choice of Einstein metric.
\end{remark}

\begin{corollary}
Let $G/K$ be a generalized flag manifold with two isotropy summands and let
$J=(+,-)$ be the invariant almost complex structure.  Then every invariant Hermitian metric \eqref{eq:metric-two} satisfies
Gray's third curvature identity.  Thus, within the invariant metric family,
\begin{equation}\label{eq:gray-whole-family1}
\{g_{\lambda}:G_3(g_{\lambda},J)=0\}
 =\{g_{\lambda}:\lambda_1>0,\ \lambda_2>0\}.
\end{equation}
    
\end{corollary}

\subsection{Nullity Kähler and distribution }

%
%
%
%
%
%
%

Let
$M=G/K$
and let $o=eK$. At the base point,
$
T_oM\simeq\got{m}.
$
For two isotropy summands,
$
\got{m}=\got{m}_1\oplus\got{m}_2,
$
since each $\got{m}_i$ is $\operatorname{Ad}(K)$-invariant, these subspaces can be transported
by the group $G$.

Define
\[
\mathcal{D}_i(gK)=(dL_g)_o(\got{m}_i).
\]
Thus
\[
TM=\mathcal{D}_1\oplus\mathcal{D}_2,
\]
and at the origin,
\[
\mathcal{D}_i(o)=\got{m}_i.
\]
Therefore, studying $\got{m}_i$ is equivalent to studying the homogeneous
distribution $\mathcal{D}_i$.

\begin{definition}[\cite{kobaiashi1}]
A distribution $\mathcal{D}\subset TM$ is integrable when, for vector fields
$X,Y\in\Gamma(\mathcal{D})$,
\[
[X,Y]\in\Gamma(\mathcal{D}).
\]
\end{definition}
By Frobenius' Theorem, this is equivalent to the local existence of submanifolds whose
tangent spaces are exactly $\mathcal{D}$.

Thus, for $\mathcal{D}_1$, we want
\[
[X,Y]\in\mathcal{D}_1,
\qquad X,Y\in\Gamma(\mathcal{D}_1).
\]
For a homogeneous distribution, this can be studied at the base point through the
$\got{m}$-component of the bracket.




\subsection*{ $\mathcal{D}_2$ is the vertical distribution}

We have
$
[\got{m}_2,\got{m}_2]\subset\mathfrak{k}.
$
Hence
$
[\got{m}_2,\got{m}_2]_{\got{m}}=0.
$
This means that, on the homogeneous space,
$
[X,Y]_o=0
$
for the invariant fields associated with $X,Y\in\got{m}_2$.
Therefore,
\[
\mathcal{D}_2\text{ is integrable}.
\]

Note that
$\mathfrak{l}=\got{k}\oplus\got{m}_2$. As $[\got{k},\got{m}_{2}]\subset \got{m}_{2}$ and $[\got{m}_{2},\got{m}_{2}]\subset \got{k}$ we obtain that $[\got{l},\got{l}]\subset \got{l}$, that is, $\got{l}$ is a subalgebra. Now consider the correspondent closed subgroup $L$ with $K\subset L \subset G$.
Thus we know that two-summand flag manifolds can often be viewed through a fibration
\[
F\longrightarrow G/K\longrightarrow G/L
\]
with fiber
$
F=L/K.
$
 See, for instance, \cite{fibracao}.
The decomposition
\[
\got{m}=\got{m}_1\oplus\got{m}_2
\]
is precisely the horizontal--vertical decomposition. In general, $\got{m}_2$ is the vertical
part.

Thus
\[
\mathcal{D}_2=\ker(d\pi)
\]
for the projection
\[
\pi:G/K\longrightarrow G/L.
\]
The leaves of $\mathcal{D}_2$ are therefore the fibers of the projection.

\subsection*{The distribution $\mathcal{D}_1$}

Here
$
[\got{m}_1,\got{m}_1]\subset\mathfrak{k}\oplus\got{m}_2.
$
For $\mathcal{D}_1$ to be integrable, we would need
$
[\got{m}_1,\got{m}_1]_{\got{m}}\subset\got{m}_1.
$
But what we have is
\[
[\got{m}_1,\got{m}_1]_{\got{m}}\subset\got{m}_2.
\]
Therefore, if there exists
\[
[X,Y]_{\got{m}}\neq0,
\qquad X,Y\in\got{m}_1,
\]
it lies in $\got{m}_2$, outside $\mathcal{D}_1$.

Thus $\mathcal{D}_1$ is not integrable whenever
$
[\got{m}_1,\got{m}_1]_{\got{m}_2}\neq0.
$
This is the typical behavior of non-symmetric two-summand flag manifolds.

We conclude that
$$
\mathcal{D}_{1} \ \mbox{is integrable} \Longleftrightarrow [\mathfrak{m}_{1},\mathfrak{m}_{1}]\subset \mathfrak{k}.
$$

\subsection*{The leaves of $\mathcal{D}_2$ are totally geodesic}

For the Levi-Civita connection of a homogeneous metric, the normal component of
$D_XY$, $X,Y\in\got{m}_2$, can be analyzed using Koszul's formula:
\[
2g(D_XY,Z)
=
g([Z,X]_{\got{m}},Y)
+g(X,[Z,Y]_{\got{m}})
+g([X,Y]_{\got{m}},Z).
\]
Take
$
X,Y\in\got{m}_2$
and
$Z\in\got{m}_1.
$
Since
$
[\got{m}_2,\got{m}_2]_{\got{m}}=0
$
and
$
[\got{m}_1,\got{m}_2]\subset\got{m}_1,
$
we have
\[
[Z,X]_{\got{m}}\in\got{m}_1,\qquad [Z,Y]_{\got{m}}\in\got{m}_1.
\]
Since $\got{m}_1\perp\got{m}_2$,
\[
g([Z,X]_{\got{m}},Y)=0,
\qquad
g(X,[Z,Y]_{\got{m}})=0.
\]
Moreover,
$
[X,Y]_{\got{m}}=0.
$
Hence
$
g(D_XY,Z)=0
$
for every $Z\in\got{m}_1$. Therefore,
\[
D_XY\in\got{m}_2.
\]
Note that the normal component of $D_{X}Y$ is
tensorial along the leaves and $G$-invariant; therefore,
it suffices to verify it at the base point.

Thus $\mathcal{D}_2$ is integrable and its leaves are totally geodesic. This holds for every invariant metric
$g$.

We can formulate the following theorem.

\begin{theorem}
Let
$
M=G/K
$
be a compact homogeneous flag manifold with two isotropy summands
$
\got{m}=\got{m}_1\oplus\got{m}_2
$
and invariant metric
$
g.
$
Then the distribution
\[
\mathcal{D}_2=G\times_K\got{m}_2
\]
is integrable and its leaves are totally geodesic.
\end{theorem}

We now relate Kähler nullity to the invariant distribution $\mathcal{D}_{2}$.   

\begin{theorem}
Let $M=G/K$ be a flag manifold with two isotropy summands and $J=(+,-)$ be the non-integrable almost complex structure. If $[\mathfrak{m}_{1},\mathfrak{m}_{1}]_{\mathfrak{m}_{2}}\ne0$, then,
$
\mathcal{K}=\left\{\begin{array}{cc}
   0,  & \lambda_{2} \neq 2\lambda_{1} \\
   \mathcal{D}_{2},  & \lambda_{2}= 2\lambda_{1}
\end{array}\right..   
$
\end{theorem}
\begin{proof}
As $[\mathfrak{m}_{1},\mathfrak{m}_{1}]_{\mathfrak{m}_{2}}\ne0$, there are roots of degree 1, $\alpha, \beta$, such that $\delta=\alpha+\beta$ is a root of degree 2 and $m_{\alpha,\beta}\ne0$. The nullity of the base point,
\[
\mathcal{K}_{o}=\{X\in\mathfrak{m} : (D_{X}J)Y=0 \text{ for all } Y\in\mathfrak{m}\},
\]
is $\operatorname{Ad}(K)$-invariant. As $\mathfrak{m}_{1}$ and $\mathfrak{m}_{2}$ are irreducible and inequivalent, each intersection $\mathcal{K}_{o}\cap\mathfrak{m}_{i} = 0$ or $\mathfrak{m}_{i}$.

For $X_{\alpha}\in\mathfrak{m}_{1}^{\mathbb{C}}$, using $\alpha+\beta-\delta=0$, we obtain $\lambda_{\alpha}=\lambda_{\beta}=\lambda_{1}$, $\lambda_{-\delta}=\lambda_{2}$ and $\epsilon_{\beta}=\epsilon_{-\delta}=+1$.

Thus
\[
-\lambda_{\alpha}+\lambda_{\beta}+\lambda_{-\delta}=\lambda_{2}\ne0
\]
and, as $m_{\alpha,\beta}\ne0$,
\[
D_{X_{\alpha}}J\ne0.
\]
Hence $\mathcal{K}_{o}\cap\mathfrak{m}_{1}=0$.

Now consider a root $\delta$ of degree 2. If $\delta+\eta+\zeta=0$ is a complementary triple, then necessarily $\deg \eta=\deg \zeta=-1$, so that $\lambda_{\eta}=\lambda_{\zeta}=\lambda_{1}$.

The metric factor is
\[
-\lambda_{2}+2\lambda_{1}.
\]
If $\lambda_{2}=2\lambda_{1}$, every potentially nonzero component
of $D_{X_{\delta}}J$ vanishes; the same argument follows for roots of degree $-2$.

Consequently, $\mathfrak{m}_{2}\subset\mathcal{K}_{o}$ and, as $\mathcal{K}_{o}\cap\mathfrak{m}_{1}=0$,
\[
\mathcal{K}_{o}=\mathfrak{m}_{2}.
\]

If $\lambda_{2}\ne2\lambda_{1}$, we use $\delta=\alpha+\beta$ again. As $\delta-\alpha=\beta$ is root, we have $m_{\delta,-\alpha}\ne0$. Since $g((D_{X_{\delta}}J)X_{-\alpha},X_{-\beta})$ is
non-zero multiple of $-\lambda_{2}+2\lambda_{1}$. Therefore $\mathfrak{m}_{2}$ is not contained in the nullity; by irreducibility,
\[
\mathcal{K}_{o}\cap\mathfrak{m}_{2}=0.
\]
We conclude that $\mathcal{K}_{o}=0$.
\end{proof}

\begin{corollary}
Since \((g,J)\) is quasi-Kähler by Lemma \ref{QK} and satisfies Gray's second identity by Theorem \ref{prop:all-gray}, Gray's Theorem implies that the Kähler nullity distribution is integrable and its leaves are Kähler.
\end{corollary}

\subsection{Examples}

\subsubsection{Flag manifold $G_{2}/U(2)$}
Let us consider the 10-dimensional flag manifold $G_{2}/U(2)$, where the isotropy group  $U(2)$ is represented by the short root. The isotropy representation of this homogeneous space decomposes into two irreducible components:
$\got{m}_{1}$ and $\got{m}_{2}$, described explicitly as 
$$
\begin{array}{l}
\got{m}_{1}=R^{+}(\alpha_{1},1)=\{\alpha_{1},\alpha_{1}+\alpha_{2},\alpha_{1}+2\alpha_{2},\alpha_{1}+3\alpha_{2}\},\\
\got{m}_{2}=R^{+}(\alpha_{1},2)=\{2\alpha_{1}+3\alpha_{2}\}.
\end{array}
$$




We will denote an invariant metric $g$ by a pair of positive numbers $(\lambda_{1},\lambda_{2})$.

We consider the invariant almost complex structure $J=J_{2}=(+,-)$.

The triples of roots which the sum is zero are

\begin{equation}
    \alpha_{1}+2\alpha_{2}+(-2\alpha_{1}-3\alpha_{2})+(\alpha_{1}+\alpha_{2})=0
\end{equation}

\begin{equation}\label{tripla1}    
\alpha_{1}+3\alpha_{2}+(-2\alpha_{1}-3\alpha_{2})+(\alpha_{1})=0.
\end{equation}
We take $\mu=2\alpha_{1}+3\alpha_{2}$. We now verify the defining condition for the
Kähler nullity
$$
X\in \mathcal{K} \Leftrightarrow D\omega(X,\cdot,\cdot)=0.
$$

Note that from triple \ref{tripla1}, $\varepsilon_{\mu-\alpha_{1}}=+1$, $\varepsilon_{-\mu}=+1$,
$\lambda_{\alpha_{1}}=\lambda_{\alpha_{1}+3\alpha_{2}}$ and $\lambda_{\mu}=\lambda_{2}$.

If $\gamma=\alpha_{1}+3\alpha_{2}$, then
\begin{eqnarray*}
&& D\omega(X_{\alpha_{1}},X_{-\mu},X_{\gamma})\\
&=&-\dfrac{im_{\alpha_{1},-\mu}}{2}(\varepsilon_{\gamma}+\varepsilon_{-\mu})(-\lambda_{\alpha_{1}}+\lambda_{\mu}+\lambda_{\gamma})\\
&=&-im_{\alpha_{1},-\mu}\lambda_{2}\neq 0.
\end{eqnarray*}

In the same way we can prove that $D\omega(X_{\alpha_{1}+\alpha_{2}},\cdot,\cdot)\neq 0$, \\ $D\omega(X_{\alpha_{1}+2\alpha_{2}},\cdot,\cdot)\neq 0$ and $D\omega(X_{\alpha_{1}+3\alpha_{2}},\cdot,\cdot)\neq 0$.

Hence $\mathcal{K}\cap \got{m}_{1}=\{0\}$.

We now determine $D_{X}J$ for $X\in \got{m}_{2}$.

Let
\[
\mu=2\alpha_1+3\alpha_2
\]
be the complementary positive root of
degree two. Every complementary root
triple involving $\mu$ can be written as
\[
\mu-\beta-\gamma=0,
\]
where
\[
\beta,\gamma\in\Pi_{M,1}^{+},
\qquad
\gamma=\mu-\beta,
\]
and $\Pi_{M,1}^{+}$ denotes the set of
positive complementary roots of degree
one. We consider only those roots
$\beta$ for which $\gamma$ is also a root.

For each such complementary root triple,
the covariant derivative of the Kähler
form is given by
\begin{eqnarray*}
&& D\omega(X_{\mu},X_{-\beta},X_{-\gamma})\\
&=&-\dfrac{im_{\mu,-\beta}}{2}(\varepsilon_{-\beta}+\varepsilon_{-\gamma})(-\lambda_{\mu}+\lambda_{\beta}+\lambda_{\gamma})\\
&=&im_{\mu,-\beta}(-\lambda_{2}+2\lambda_{1})=0 \Leftrightarrow \lambda_{2}=2\lambda_{1}.
\end{eqnarray*}

In the same way, we conclude that $D\omega(X_{-(2\alpha_{1}+3\alpha_{2})},\cdot,\cdot)=0 \Leftrightarrow \lambda_{2}=2\lambda_{1}$.

Hence, $\mathfrak{m}_{2}\subset \mathcal{K}$ when $\lambda_{2}=2\lambda_{1}$.
 
In summary, the Kähler nullity of $G_{2}/U(2)$ is $\mathcal{K}=\left\{\begin{array}{ccl}
0 & \mbox{for} & \lambda_{2}\neq 2\lambda_{1}\\
\got{m}_{2} & \mbox{for} & \lambda_{2}=2\lambda_{1}
\end{array}\right.
.$

Geometrically, $\mathcal{K}=\mathfrak{m}_{2}$ is identified with the vertical distribution of the twistor fibration
$G_{2}/U(2) \longrightarrow G_{2}/SO(4)$.


Since $(G_2/U(2),g,J_2)$ is quasi-K\"ahler by
Lemma~3.1 and satisfies Gray's second curvature
identity by Theorem~3.2, Theorem~2.12 yields
\[
\text{quasi-K\"ahler}+\text{Gray II}
\Longrightarrow
\mathcal{K}\text{ is integrable}
\Longrightarrow
\text{the leaves are K\"ahler}.
\]

\subsubsection{The flag manifold $\mathbb{CP}^{3}=Sp(2)/(Sp(1)\times U(1))$}\label{section-CP3}

Let us consider the six-dimensional homogeneous space
\[
\mathbb{CP}^{3}
=Sp(2)/(Sp(1)\times U(1)).
\]
We work with the complexified Lie algebra
\(\mathfrak{sp}(4,\mathbb C)\) of the compact real
Lie algebra \(\mathfrak{sp}(2)\). Its standard
Cartan subalgebra is
\[
\mathfrak h=
\left\{
\operatorname{diag}(x,y,-x,-y):
x,y\in\mathbb C
\right\}.
\]
Equivalently, its elements can be written as
\[
H=\operatorname{diag}(\Lambda,-\Lambda),
\qquad
\Lambda=\operatorname{diag}(x,y).
\]

One can describe precisely each irreducible isotropy summand $\mathfrak{m}_1, \mathfrak{m}_2$ in terms of the root space decomposition of the Lie algebra $\mathfrak{sp}(2)$.  We will use the following basis of $\mathfrak{sp}(2)$: let $E_{ij}$ $(1 \leq i,j \leq 4)$  be the $4 \times 4$ matrix with $1$ in the $ij$-position and zero otherwise. We define 
$$
X_{12}  =  \dfrac{\sqrt{3}}{6}(E_{12}-E_{43}), \,\,\, X_{-12}  =  \dfrac{\sqrt{3}}{6}(E_{21} - E_{34}), \,\,\, X_{12}^{+}  =  \dfrac{\sqrt{3}}{6}(E_{14}+E_{23}), 
$$
$$
X_{-12}^{+}  =  \dfrac{\sqrt{3}}{6}(E_{32}+E_{41}), \,\,\, X_{11}  =  \dfrac{\sqrt{6}}{6}(E_{13}),  \,\,\, X_{-11}  =  \dfrac{\sqrt{6}}{6}(E_{31}).
$$
%


The matrices defined above are root vectors
in \(\mathfrak{sp}(4,\mathbb C)\), associated
with the Cartan subalgebra \(\mathfrak h\)
specified above.

Let us consider the linear functional given by 
$$
e_{1}:\Lambda=Diag\{x,y\}\longrightarrow x \ \ and \ \ 
e_{2}:\Lambda=Diag\{x,y\}\longrightarrow y.
$$

We have the following positive complementary roots for $\mathbb{CP}^{3}$: 
$$
\alpha_{1}=e_{1}-e_{2}, \ \ \alpha_{2}=e_{1}+e_{2} \ \  \mbox{e} \ \ \alpha_{3}=2e_{1}.
$$


The isotropy representation decomposes in the following way
\begin{eqnarray*}
\mathfrak{m}&=&\mathfrak{m}_{1}\oplus \mathfrak{m}_{2}\\
&=&(\mathfrak{u}_{\alpha_{1}}\oplus \mathfrak{u}_{\alpha_{2}}) \oplus  \mathfrak{u}_{\alpha_{3}}.
\end{eqnarray*}

We consider the invariant almost complex structure $J=J_{2}=(+,-)$.

Note that
\[
[X_{\alpha_1},X_{\alpha_2}]=
m_{\alpha_{1},\alpha_{2}} X_{\alpha_3}.
\]
Since
$
X_{\alpha_1},X_{\alpha_2}\in\got{m}_1,
$
and
$
X_{\alpha_3}\in\got{m}_2,
$
we have
$
[\mathcal{D}_1,\mathcal{D}_1]\not\subset\mathcal{D}_1.
$
Therefore $\mathcal{D}_1$ is not integrable.

On the other hand,
$
[\got{m}_2,\got{m}_2]\subset\mathfrak{k},
$
and therefore $\mathcal{D}_2$ is integrable.


For $\mathcal{D}_2$, Koszul's formula gives, for
$X,Y\in\got{m}_2$ and $Z\in\got{m}_1$,
\begin{equation}\label{eq:distribution}
g(D_XY,Z)=0.
\end{equation}
Indeed, $[X_{\alpha_{3}},X_{\alpha_{3}}]_{\mathfrak{m}}=[X_{\alpha_{3}},X_{\alpha_{-3}}]_{\mathfrak{m}}=0$, $[X_{\alpha_{1}},X_{\alpha_{-3}}]_{\mathfrak{m}}\subset\got{m}_1$, $[X_{\alpha_{2}},X_{-\alpha_{-3}}]_{\mathfrak{m}}\subset\got{m}_1$  and
\[
[\got{m}_2,\got{m}_2]_{\got{m}}=0,
\qquad
[\got{m}_1,\got{m}_2]\subset\got{m}_1. 
\]
As
$
\got{m}_1\perp\got{m}_2
$, the equation \ref{eq:distribution} follows.
Consequently,
\[
D_XY\in\got{m}_2.
\]
Thus $\mathcal{D}_2$ is integrable and totally geodesic. 

There is only a 
triple of roots which the sum is zero:
$$
\alpha_{1}+\alpha_{2}-\alpha_{3}=0.
$$

We observe that
\begin{eqnarray*}
&& D\omega(X_{1},X_{2},X_{-3})\\
&=&-\dfrac{im_{1,2}}{2}(\varepsilon_{2}+\varepsilon_{-3})(-\lambda_{1}+\lambda_{1}+\lambda_{2})\\
&=&-im_{1,2}(\lambda_{2}) \neq 0.
\end{eqnarray*}
Hence $D\omega(X_{1},\cdot,\cdot)\neq 0$

In the same way $D\omega(X_{2},\cdot,\cdot)\neq 0$, $D\omega(X_{-2},\cdot,\cdot)\neq 0$ and $D\omega(X_{-1},\cdot,\cdot)\neq 0$.
Moreover,
\begin{eqnarray*}
&& D\omega(X_{-3},X_{1},X_{2})\\
&=&-\dfrac{im_{-3,1}}{2}(\varepsilon_{1}+\varepsilon_{2})(-\lambda_{2}+\lambda_{1}+\lambda_{1})\\
&=&-im_{-3,1}(-\lambda_{2}+2\lambda_{1})=0\Longleftrightarrow\lambda_{2}=2\lambda_{1}.
\end{eqnarray*}
and
$D\omega(X_{3},\cdot,\cdot)= 0 \Longleftrightarrow \lambda_{2}=2\lambda_{1}$.

Hence, $\mathfrak{m}_{2}\subset \mathcal{K}$ when $\lambda_{2}=2\lambda_{1}$.

Therefore, the Kähler nullity of $\mathbb{CP}^{3}$ is $\mathcal{K}=\left\{\begin{array}{ccl}
0 & \mbox{for} & \lambda_{2}\neq 2\lambda_{1}\\
\mathcal{D}_2 & \mbox{for} & \lambda_{2}=2\lambda_{1}
\end{array}\right..
$

The leaves of $\mathcal{D}_{2}$ are the fibers of the fibration
$$
\mathbb{CP}^{1}\longrightarrow \mathbb{CP}^{3}  \longrightarrow S^{4}.
$$
As in the example for $G_{2}/U(2)$, the leaves are $\mathbb{CP}^{1}$.

\subsection{The integrable case}\label{integrable}

We now briefly discuss the integrable invariant complex structure
\[
J_{1}=J_+=(+,+).
\]
Our purpose is not to repeat the arguments developed for the non-integrable structure ($J=(+,-)$), but rather to point out the differences between the two cases.

For $(J_+)$, the signs associated with complementary roots of degrees $(1,-1,2,-2)$ are
\[
\begin{array}{c|cccc}
\deg\alpha & 1 & -1 & 2 & -2\\ \hline
\varepsilon_\alpha & +1 & -1 & +1 & -1.
\end{array}
\]

We now verify Gray's second curvature
identity for the integrable invariant
complex structure $J_+$.

As in the proof of Theorem~3.2, the
isotropy-torus invariance of the curvature
tensor implies that a potentially nonzero
component must involve complementary
root degrees whose sum is zero.

The possible degree multisets are
precisely those listed in~(9):
\[
(1,1,-1,-1),\quad
(2,2,-2,-2),\quad
(1,-1,2,-2).
\]
For $J_+$, each multiset contains exactly
two positive and two negative $J$-signs.

Consequently, the coefficient
\[
1+\varepsilon_\alpha\varepsilon_\beta
+\varepsilon_\alpha\varepsilon_\gamma
+\varepsilon_\alpha\varepsilon_\delta
\]
in the Gray-II defect vanishes for every
potentially nonzero curvature component.
Hence
\[
G_2(g_\lambda,J_+)=0
\]
for every $\lambda_1,\lambda_2>0$.

The set such that $g$ satisfies the second Gray identity is  $$\{
\lambda_1Q|{\mathfrak m_1}
+
\lambda_2Q|{\mathfrak m_2},
\qquad \lambda_1,\lambda_2>0.
\}
$$

Thus Gray's second identity does not distinguish the integrable and non-integrable invariant almost complex structures within the family of invariant metrics.

The quasi-K\"ahler condition, however, behaves
differently. Since
$[\mathfrak{m}_1,\mathfrak{m}_1]_{\mathfrak{m}_2}
\neq 0$, there exists a complementary root triple
with nonzero structure constant satisfying

\[
\alpha+\beta-\delta=0,
\]
with
\[
\deg\alpha=\deg\beta=1,
\qquad
\deg\delta=2.
\]
For $(J_+)$, we have
\[
\varepsilon_\alpha=\varepsilon_\beta=1,
\qquad
\varepsilon_{-\delta}=-1.
\]
Hence the quasi-Kähler condition becomes
\[
\lambda_1+\lambda_1-\lambda_2=0,
\]
and therefore
\[
\lambda_2=2\lambda_1.
\]
Since $(J_+)$ is integrable, this is precisely the Kähler condition. Thus
\[
(g_\lambda,J_+)\text{ is quasi-Kähler}
\quad\Longleftrightarrow\quad
(g_\lambda,J_+)\text{ is Kähler}
\quad\Longleftrightarrow\quad
\lambda_2=2\lambda_1.
\]

We now turn to the Kähler nullity. Recall that, for a complementary root triple
\[
\alpha+\beta+\gamma=0,
\]
one has, up to a nonzero normalization factor,
\[
g((D_{X_\alpha}J)X_\beta,X_\gamma)
=m_{\alpha,\beta}
(\varepsilon_\beta+\varepsilon_\gamma)
(-\lambda_\alpha+\lambda_\beta+\lambda_\gamma).
\]
If $X_{\alpha}\in \got{m}_{1}$, then 
$\varepsilon_\beta+\varepsilon_\gamma=1-1=0$ because the triples have degrees $1+1-2=0$ or $1-2+1=0$
and therefore
\[
(D_{X_\alpha}J_+)=0.
\]
Hence
\[
\mathfrak m_1\subset\mathcal K.
\]

On the other hand, if $(X_\delta\in\mathfrak m_2)$, the relevant triples have degrees
$
2-1-1=0.
$
Now
$
\varepsilon_{-\alpha}+\varepsilon_{-\beta}=-2,
$
while the remaining metric factor is
$
-\lambda_2+2\lambda_1.
$
Therefore
\[
(D_{X_\delta}J_+)=0
\quad\Longleftrightarrow\quad
\lambda_2=2\lambda_1.
\]

We obtain
\[
\mathcal K=
\begin{cases}
\mathfrak m_1, & \lambda_2\neq2\lambda_1,\\
\mathfrak m_1\oplus\mathfrak m_2=\mathfrak m,
& \lambda_2=2\lambda_1.
\end{cases}
\]

Thus the integrable and non-integrable structures exhibit complementary behavior. For ($J=(+,-)$), the Kähler nullity is trivial away from the distinguished ratio and equals $(\mathfrak m_2)$ when ($\lambda_2=2\lambda_1)$. For ($J_+=(+,+)$), the nullity equals ($\mathfrak m_1$) away from that ratio and becomes the whole tangent space at the Kähler metric.

Moreover, whenever
\[
[\mathfrak m_1,\mathfrak m_1]_{\mathfrak m_2}\neq0,
\]
the distribution determined by ($\mathfrak m_1$) is not integrable. Hence, for
\[
\lambda_2\neq2\lambda_1,
\]
the structure ($(g_\lambda,J_+)$) provides an example satisfying Gray's second identity whose Kähler nullity is not integrable. This does not contradict Gray's integrability Theorem, since the quasi-Kähler hypothesis fails in this range.

%
%
%


We can summarize these facts in the table \ref{tab:comparison_gray_nullity}.

\begin{table}[htbp]
\centering
\caption{Comparison of quasi-Kählerity, Gray II identity, and Kähler nullity for invariant almost Hermitian structures on generalized flag manifolds with two isotropy summands.}
\label{tab:comparison_gray_nullity}
\begin{tabular}{lccc}
\toprule
\textbf{Structure} & \textbf{quasi-Kähler?} & \textbf{Gray II} & \textbf{Kähler nullity} \\
\midrule
$J_2 = (+,-), \quad \lambda_2 \neq 2\lambda_1$ & yes & yes & $\{0\}$ \\
[1.5mm]
$J_2 = (+,-), \quad \lambda_2 = 2\lambda_1$ & yes & yes & $\mathcal{D}_2$ \\
[1.5mm]
$J_+ = (+,+), \quad \lambda_2 \neq 2\lambda_1$ & no & yes & $\mathcal{D}_1$ \\
[1.5mm]
$J_+ = (+,+), \quad \lambda_2 = 2\lambda_1$ & Kähler & yes & $TM$ \\
\bottomrule
\end{tabular}
\end{table}
Hence $G_2$ alone does not imply integrability of the K\"ahler nullity.

\section{ Second Gray identity along homogeneous Ricci flow}
\label{sec:application}

In this section, we investigate the behavior of invariant almost Hermitian structures satisfying Gray's second curvature identity under the homogeneous Ricci flow.

Recall that an invariant metric on a generalized flag manifold with two isotropy summands is given by
$$
g_\lambda=\lambda_1Q|{\mathfrak m_1}+\lambda_2Q|{\mathfrak m_2},
\qquad \lambda_1,\lambda_2>0.
$$

\subsection{Ricci flow and the ratio of metric parameters}
\label{subsec:ratio-flow}

The invariant Ricci flow preserves the space of metrics \eqref{eq:metric-two}. This finite-dimensional reduction is consistent with the general framework of homogeneous bracket flow developed by Lauret \cite{Lauret2013,Lauret2016}. The same holds for the normalized Ricci flow, since its additional scalar term rescales both summands equally.
 The homogeneous Ricci flow on generalized flag manifolds was studied in \cite{AnastassiouChrysikos2011}; see also \cite{buzano,GramaMartins2012}. We will use the framework provided in \cite{AnastassiouChrysikos2011}. The Einstein metrics in this family were determined in \cite{AC2011}. Put
\begin{equation}\label{eq:ratio}
u=\frac{\lambda_2}{\lambda_1},
\qquad d_i=\dim\m_i.
\end{equation}
In the normalization $Q=-B$, put
\begin{equation}\label{eq:A}
A
=\frac{d_1d_2}{d_1+4d_2}.
\end{equation}

By ~\cite[Proposition~2(1)]{AnastassiouChrysikos2011},
the eigenvalues $\rho_i$ of the Ricci operator, defined by
$\operatorname{Ric}^{\mathrm{op}}|_{\mathfrak m_i}
= \rho_i\,\operatorname{Id}_{\mathfrak m_i}$, are given by

\begin{align}
\rho_1&=\frac{1}{2\lambda_1}
             -\frac{A\lambda_2}{2d_1\lambda_1^2},
 \label{eq:rho1}\\
\rho_2&=\frac{1}{2\lambda_2}
  +\frac{A}{4d_2}\left(\frac{\lambda_2}{\lambda_1^2}
                               -\frac{2}{\lambda_2}\right).
 \label{eq:rho2}
\end{align}

The constant \(A\) in (15) is the structural
constant usually denoted by \([112]\), computed
with respect to \(Q=-B\). More precisely, if
\(\{e_\alpha^i\}\) is a \(Q\)-orthonormal basis
of \(\mathfrak m_i\), then
\[
[112]=
\sum_{\alpha,\beta,\gamma}
Q\bigl([e_\alpha^1,e_\beta^1]_{\mathfrak m_2},
e_\gamma^2\bigr)^2.
\]
The bracket relations
\[
[\mathfrak m_1,\mathfrak m_1]
\subset\mathfrak k\oplus\mathfrak m_2,
\qquad
[\mathfrak m_1,\mathfrak m_2]
\subset\mathfrak m_1,
\qquad
[\mathfrak m_2,\mathfrak m_2]
\subset\mathfrak k
\]
imply that this is the only nonzero structural
constant involving three isotropy summands,
up to permutation.

Applying the homogeneous Ricci curvature formula
in~\cite[Proposition~1(1)]{AnastassiouChrysikos2011},
we obtain

\[
\operatorname{Ric}|_{\mathfrak m_1}
=
\left(
\frac12-\frac{A\lambda_2}{2d_1\lambda_1}
\right)Q|_{\mathfrak m_1}
\]
and
\[
\operatorname{Ric}|_{\mathfrak m_2}
=
\left(
\frac12-\frac{A}{2d_2}
+\frac{A\lambda_2^2}{4d_2\lambda_1^2}
\right)Q|_{\mathfrak m_2}.
\]
Dividing these coefficients by
\(\lambda_1\) and \(\lambda_2\),
respectively, gives the Ricci-operator
eigenvalues (16)--(17).

Here it is essential to distinguish the Ricci tensor coefficients from the
operator eigenvalues: in the notation $\operatorname{Ric}|_{\m_i}=r_iQ|_{\m_i}$,
one has $\rho_i=r_i/\lambda_i$.

We emphasize the conventions used in the
dynamical analysis below. The Ricci flow
is written as
\[
\frac{\partial g}{\partial t}
=-2\operatorname{Ric}(g),
\]
and \(\rho_i\) denotes the eigenvalue of
the Ricci endomorphism, rather than the
coefficient of the Ricci tensor relative
to \(Q\). Thus the unnormalized parameter
equations are
\[
\frac{d\lambda_i}{dt}
=-2\lambda_i\rho_i.
\]
The ratio variable \(u=\lambda_2/\lambda_1\)
eliminates the common scaling direction,
and the time change
\[
\frac{d\tau}{dt}=\frac1{\lambda_1(t)}
\]
produces the autonomous equation (19).
These conventions are used throughout
our comparison with the homogeneous
Ricci-flow framework of \cite{GramaMartins2012}.

For a normalized Ricci flow, write
$\dot\lambda_i=-2\lambda_i\rho_i+c(t)\lambda_i$, with the same function
$c(t)$ on both summands.  Differentiating \eqref{eq:ratio}, the normalization
cancels and gives
\begin{equation}\label{eq:ratio-derivative}
\dot u
=\frac{\dot\lambda_2\lambda_1-\lambda_2\dot\lambda_1}
       {\lambda_1^2}
=2u(\rho_1-\rho_2).
\end{equation}
Substitute $\lambda_2=u\lambda_1$ into \eqref{eq:rho1}--\eqref{eq:rho2}:
\begin{align*}
\lambda_1(\rho_1-\rho_2)
&=\frac12-\frac{Au}{2d_1}-\frac{1}{2u}
   -\frac{A}{4d_2}\left(u-\frac2u\right)\\
&=-\frac{(u-2)((d_1+2d_2)u-4d_2)}
            {4u(d_1+4d_2)}.
\end{align*}
Use the increasing time parameter $d\tau/dt=1/\lambda_1(t)$.  Equation
\eqref{eq:ratio-derivative} becomes
\begin{equation}\label{eq:general-flow}
\frac{du}{d\tau}
=-\frac{(u-2)((d_1+2d_2)u-4d_2)}{2(d_1+4d_2)}.
\end{equation}
Thus the two invariant Einstein ratios are
\begin{equation}\label{eq:einstein-ratios}
u_E^{(K)}=2,\qquad
u_E^{(N)}=\frac{4d_2}{d_1+2d_2}<2.
\end{equation}
The superscript $(K)$ refers to the K\"ahler--Einstein metric for the
integrable complex structure; $(N)$ denotes the other invariant Einstein
metric.  Neither superscript characterizes the Gray II locus of
Theorem~\ref{prop:all-gray}.

The right-hand side of \eqref{eq:general-flow} is negative for
$0<u<u_E^{(N)}$, positive for $u_E^{(N)}<u<2$, and negative for $u>2$.
Therefore $u_E^{(N)}$ is repelling and $u=2$ is attracting in the
projectivized positive metric space.  The statement concerns the ratio of
parameters, not convergence of the unnormalized metric without rescaling.

\begin{proposition}
\label{prop:ratio-trajectories}
Put
\[
a=\frac{4d_2}{d_1+2d_2},\qquad
c=\frac{d_1+2d_2}{2(d_1+4d_2)}.
\]
Then $0<a<2$, $c>0$, and the reduced equation is
$u'=-c(u-a)(u-2)$.  If $u(0)>a$, the solution exists for all
positive rescaled time and $u(\tau)\to2$.  If $0<u(0)<a$, the ratio
decreases to the boundary $u=0$ in finite rescaled time.  The constant
solutions are $u=a$ and $u=2$.
\end{proposition}

\begin{proof}
For a nonconstant solution, separation of variables gives
\begin{equation}\label{eq:integrated-ratio}
\frac{d}{d\tau}\log\left|\frac{u-2}{u-a}\right|
=-c(2-a),\qquad
\frac{u(\tau)-2}{u(\tau)-a}
=\frac{u(0)-2}{u(0)-a}e^{-c(2-a)\tau}.
\end{equation}
The sign of $u'$ shows that a solution starting in $(a,2)$ increases,
whereas one starting above $2$ decreases.  Both approach $2$; this also
follows from \eqref{eq:integrated-ratio}.  On $(0,a)$ the right-hand side
of the ODE is strictly negative.  It has a finite negative value at $u=0$
in its polynomial extension, so a solution in this interval reaches the
boundary in finite $\tau$.  We stop the geometric flow interpretation at
that boundary, where the metric ceases to be positive definite.
\end{proof}

In particular, Gray II imposes no extra restriction on initial data in
this portrait: each positive initial ratio is a Gray metric, including
those trajectories moving towards the boundary.  A positive rescaling of
time changes neither their ordering nor the equilibrium ratios.

For reference, under the volume normalization $\operatorname{Vol}(g(t))$
constant, the parameter ratio still follows \eqref{eq:general-flow}.  Up
to a constant depending on the chosen volume, the two metric coefficients
can be recovered from the ratio by
\begin{equation}\label{eq:recover-lambdas}
\lambda_1=C\,u^{-d_2/(d_1+d_2)},\qquad
\lambda_2=C\,u^{d_1/(d_1+d_2)}.
\end{equation}

For \(0<u(0)<a\), the solution reaches
\(u=0\) at a finite rescaled time
\(\tau_*\). Under volume normalization,
equation (22) gives
\[
\lambda_1=C u^{-p},
\qquad
p=\frac{d_2}{d_1+d_2}<1.
\]
Since \(u'(\tau)\to-2ca<0\) as
\(\tau\to\tau_*^-\), we have
\(u(\tau)\asymp\tau_*-\tau\). Hence
\[
t_*-t=\int_\tau^{\tau_*}\lambda_1(s)\,ds
<\infty,
\]
because \(p<1\). Thus the degeneration
also occurs in finite geometric time,
with \(\lambda_1\to\infty\) and
\(\lambda_2\to0\).

Indeed $\operatorname{Vol}(g)$ is proportional to
$\lambda_1^{d_1/2}\lambda_2^{d_2/2}$, and $\lambda_2=u\lambda_1$.
This explains why stationary directions of the projectivized flow give
fixed metrics after a suitable normalization, while the unnormalized
metrics may still change by scale.

\subsection{The Type I family}
\label{subsec:typeI}

For
\begin{equation}\label{eq:typeI-space}
M_{m,k}=\frac{SO(2n+1)}{U(m)\times SO(2k+1)},
\qquad n=m+k,\qquad m\geq2,\quad k\geq0,
\end{equation}
the summand dimensions are
\begin{equation}\label{eq:typeI-dimensions}
d_1=2m(2k+1),\qquad d_2=m(m-1).
\end{equation}
Hence $d_1+2d_2=2m(m+2k)$ and
$d_1+4d_2=2m(2n-1)$.  Substitution into
\eqref{eq:general-flow} gives
\begin{equation}\label{eq:typeI-flow}
\frac{du}{d\tau}
=-\frac{(u-2)((m+2k)u-2(m-1))}{2(2n-1)},
\qquad
u_E^{(N)}=\frac{2(m-1)}{m+2k},\quad u_E^{(K)}=2.
\end{equation}
All the trajectories with $u>0$ satisfy Gray II by
Theorem~\ref{prop:all-gray}.  In particular, the vector field is generally
nonzero at a Gray metric; its zeros in \eqref{eq:typeI-flow} identify Einstein
metrics, not special Gray directions.

One can read the dependence on $(m,k)$ directly from the second
equilibrium.  Increasing $k$ with $m$ fixed moves $u_E^{(N)}$ towards
zero.  Its position relative to $u=1$ is determined without solving the
ODE:
\begin{equation}\label{eq:typeI-compare-one}
u_E^{(N)}-1=\frac{m-2k-2}{m+2k}.
\end{equation}
Hence $u_E^{(N)}<1$, $=1$, or $>1$ according as $m<2k+2$,
$m=2k+2$, or $m>2k+2$.  Gray II holds on both sides of each of
these values.

For the case \(m=2\) and \(k=1\), namely
\(\mathrm{SO}(7)/(\mathrm{U}(2)\times\mathrm{SO}(3))\),
we have \(d_1=12\) and \(d_2=2\). Thus
\begin{equation}\label{eq:SO7-flow}
\frac{du}{d\tau}=-\frac{(u-2)(2u-1)}5,
\qquad u_E^{(N)}=\frac12,\quad u_E^{(K)}=2.
\end{equation}
The direction $u=1/2$ is repelling; $u=2$ attracts ratios from either side.
Both, as well as every ratio between or beyond them, satisfy Gray II.

For this example, the complete phase line
of the projectivized ratio dynamics is

\begin{center}
\begin{tikzpicture}[>=stealth,scale=0.95]
\draw[->,thick] (0,0)--(10.2,0) node[right] {$u>0$};

\fill (2.4,0) circle (2pt);
\fill (7.3,0) circle (2pt);

\node[below] at (2.4,0) {$\frac12$};
\node[below] at (7.3,0) {$2$};

\draw[->,thick] (1.7,0.3)--(0.65,0.3);

\draw[->,thick] (3.0,0.3)--(4.5,0.3);
\draw[->,thick] (5.4,0.3)--(6.8,0.3);

\draw[->,thick] (9.1,0.3)--(7.8,0.3);

\node[above] at (2.4,0.5) {Einstein, repelling};
\node[above] at (7.3,0.5) {Einstein, attracting};

\node[below] at (5.1,-0.6) {Gray II holds along the entire positive axis};
\end{tikzpicture}
\end{center}

The horizontal coordinate is the ratio $u$, not the distance in a metric
on the space of brackets.  The arrows encode the sign of $du/d\tau$.

\subsection{The flag manifold $\mathbb{CP}^{3}=Sp(2)/(Sp(1)\times U(1))$}
\label{subsec:sp2}

We retain the root-space decomposition and
the labeling of the isotropy summands
introduced in Section \ref{section-CP3}.

With the numbering $\m_1$ of real dimension $4$ and $\m_2$ of real
dimension $2$ used above, $d_1=4$, $d_2=2$, and $A=2/3$.
To make the labeling unambiguous, take simple roots
$\alpha=e_1-e_2$ and $\beta=2e_2$ for the root system of type $C_2$.
Here $\beta$ is an \emph{isotropy root}.  The positive complementary
roots are $\alpha$, $\alpha+\beta=e_1+e_2$, and
$2\alpha+\beta=2e_1$.  The first two generate $\m_1$, and the last
generates $\m_2$.  Thus $\dim\m_1=4$ and $\dim\m_2=2$, as required;
counting $2e_2$ as a root of $\m_2$ would change the calculation.
Equations \eqref{eq:rho1}--\eqref{eq:rho2} give, in the \emph{same}
normalization $Q=-B$,
\begin{align}
\rho_1&=\frac{1}{2\lambda_1}
             -\frac{\lambda_2}{12\lambda_1^2},
 \label{eq:sp2-rho1}\\
\rho_2&=\frac{1}{3\lambda_2}
             +\frac{\lambda_2}{12\lambda_1^2}.
 \label{eq:sp2-rho2}
\end{align}
In particular,
\begin{equation}\label{eq:sp2-rho-difference}
\rho_1-\rho_2
=\frac{1}{\lambda_1}\left(\frac12-\frac{1}{3u}-\frac{u}{6}\right)
=-\frac{(u-1)(u-2)}{6\lambda_1u}.
\end{equation}
The reduced equation and its Einstein ratios are consequently
\begin{equation}\label{eq:sp2-flow}
\frac{du}{d\tau}=-\frac{(u-1)(u-2)}{3},
\qquad u_E^{(N)}=1,\quad u_E^{(K)}=2.
\end{equation}
The arrows point towards decreasing $u$ for $0<u<1$, towards increasing
$u$ for $1<u<2$, and towards decreasing $u$ for $u>2$.

The nearly K\"ahler construction of \cite[Proposition~3.1]{Nagy2002}
reverses the complex structure on the vertical summand of a K\"ahler
fibration and halves the metric on that summand.  Here the vertical
summand is $\m_2$; starting from the K\"ahler metric with ratio $2$
gives a nearly K\"ahler metric with ratio $1$.  The coincidence of this
ratio with $u_E^{(N)}$ is consistent with the six-dimensional homogeneous
example.  Nearly K\"ahler is a more restrictive property than Gray II in
this family.

In this example the distinction has a particularly simple expression:
the nearly K\"ahler Einstein ratio $u=1$ is a repelling stationary
direction of the normalized ratio dynamics, while the K\"ahler--Einstein
ratio $u=2$ is attracting.  Initial values in $(1,2)$ move from near the
first direction towards the second.  At every intervening time, the
almost complex structure under discussion remains the same
\emph{nonintegrable} $J=(+,-)$ and Gray II remains valid.  The label
``K\"ahler--Einstein ratio'' at $u=2$ refers to a \emph{different},
integrable invariant complex structure and does not turn this fixed $J$
into a K\"ahler structure.

\begin{remark}\label{subsec:bracket-conclusion}
We now check directly why the ratio equation is also the equation seen in
the bracket representation.  Let $\mu_{11}^{2}$ denote the part of the
homogeneous bracket with two inputs in $\m_1$ and output in $\m_2$;
similarly, let $\mu_{12}^{1}$ have inputs in $\m_1,\m_2$ and output in
$\m_1$.  When a metric $g_{\lambda}$ is expressed in an orthonormal
frame for the fixed background form $Q$, their respective scalar
\emph{metric factors} are
\begin{equation}\label{eq:bracket-factors}
x=\frac{\sqrt{\lambda_2}}{\lambda_1},\qquad
y=\frac{1}{\sqrt{\lambda_2}},\qquad
\frac{x}{y}=\frac{\lambda_2}{\lambda_1}=u.
\end{equation}
The factors multiply the fixed algebraic structure constants of the two
bracket components.  Other components, including those with output in
$\kalg$, remain part of the full homogeneous bracket; the two displayed
factors suffice to track the projective metric parameter.

Indeed, if the bracket equation is written
$\dot\mu=-\pi(\Ricop)\mu$, a diagonal operator with eigenvalues
$\rho_1,\rho_2$ acts on $\mu_{11}^{2}$ with weight
$\rho_2-2\rho_1$ and on $\mu_{12}^{1}$ with weight
$\rho_1-(\rho_1+\rho_2)=-\rho_2$.  Consequently, for unnormalized
bracket flow,
\begin{equation}\label{eq:bracket-factor-flow}
\frac{\dot x}{x}=2\rho_1-\rho_2,
\qquad
\frac{\dot y}{y}=\rho_2,
\qquad
\frac{d}{dt}\log\frac{x}{y}=2(\rho_1-\rho_2).
\end{equation}
This is exactly \eqref{eq:ratio-derivative}.  If the normalized metric flow is written as
\(\dot\lambda_i=-2\lambda_i\rho_i+c(t)\lambda_i\), then the
corresponding bracket-flow operator is
\(\operatorname{Ric}_{\mathrm{op}}-\frac{c(t)}{2}\operatorname{Id}\).
Consequently, both logarithmic derivatives in \eqref{eq:bracket-factor-flow} acquire the
same additional term \(-c(t)/2\), which cancels in
\(\frac{d}{dt}\log(x/y)\). Hence the same reduced
equation governs normalized bracket flow.  The identity $x/y=u$ is the
explicit change of variables needed to compare the two descriptions;
using $x$ or $y$ alone as if it were $u$ would give an incorrect scalar
ODE.

The Ricci operator is block scalar and commutes with $J$.  This fact alone
would not prove preservation of a curvature identity under an arbitrary
bracket flow: the Ricci operator is symmetric, whereas infinitesimal unitary
changes of frame are skew-symmetric.  Here preservation has a direct reason.
The homogeneous Ricci flow stays within the invariant metric family
\eqref{eq:metric-two}, and Theorem~\ref{prop:all-gray} shows that the
Gray defect vanishes identically on that family.  Therefore along every such
trajectory,
\begin{equation}\label{eq:gray-along-flow}
G_2(g(t),J)=0,
\qquad \frac{d}{dt}G_2(g(t),J)=0.
\end{equation}
Under the standard homogeneous bracket-flow correspondence
\cite{Lauret2013,Lauret2016}, the compatible bracket and almost complex
structure represent the same invariant almost Hermitian geometry up to the
corresponding equivariant changes of frame.  Thus \eqref{eq:gray-along-flow}
also holds for the Gray defect of the represented bracket-flow solution.
This is a statement about the invariant two-parameter family and its induced
bracket-flow realization; no tangency assertion for arbitrary brackets is
required.

\end{remark}

\begin{remark}
\label{rem:Kähler-nullity-ricci-flow}

The preceding results reveal a geometric phenomenon
concerning the behavior of Kähler nullity along
homogeneous Ricci flow.

Let $J_{2}= (+,-)$ be the fixed non-integrable invariant
almost complex structure, and assume that
\[
[\mathfrak m_1,\mathfrak m_1]_{\mathfrak m_2}
\neq 0.
\]
By Theorem ~3.7, the Kähler nullity of
$(g_\lambda,J_{2})$ satisfies
\[
\mathcal K(g_\lambda,J_{2})=
\begin{cases}
\{0\}, & \lambda_2\neq 2\lambda_1,\\[2mm]
D_2, & \lambda_2=2\lambda_1.
\end{cases}
\]

On the other hand, the projectivized homogeneous
Ricci flow is governed by
\[
\frac{du}{d\tau}
=-c(u-a)(u-2),
\qquad
u=\frac{\lambda_2}{\lambda_1},
\]
where
\[
a=\frac{4d_2}{d_1+2d_2},
\qquad
c=\frac{d_1+2d_2}{2(d_1+4d_2)}>0.
\]

Since $u=2$ is an equilibrium, uniqueness of
solutions implies that a nonstationary trajectory
cannot reach this value in finite time.
Consequently, if $u(0)\neq 2$, then
\[
\mathcal K(g(\tau),J_{2})=\{0\}
\]
throughout its finite interval of existence.

Nevertheless, whenever $u(0)>a$ and $u(0)\neq 2$,
the projectivized trajectory satisfies
\[
\lim_{\tau\to\infty}u(\tau)=2.
\]
Under volume normalization, the corresponding
metrics converge to an invariant Einstein metric
$g_\infty$ with ratio $u=2$.

Therefore,
\[
\lim_{\tau\to\infty}
\dim\mathcal K(g(\tau),J_{2})
=0,
\]
whereas
\[
\dim\mathcal K(g_\infty,J_{2})=d_2>0.
\]

Thus, although the Kähler nullity remains trivial
at every finite time, its dimension increases
at the limiting metric.

Notice that $g_\infty$ is Kähler--Einstein with
respect to the integrable invariant complex
structure $J_+=(+,+)$, but the fixed structure
$J_{2}$ remains non-integrable. The appearance
of the nontrivial nullity distribution $D_2$
therefore does not represent convergence to
a Kähler structure with respect to $J_{2}$.

This illustrates that preservation of Gray's
second curvature identity under homogeneous
Ricci flow does not imply constancy of the
Kähler nullity dimension under convergence
of invariant metrics. In particular, the
dimension of the nullity may increase in
the asymptotic limit without changing at
any finite time.
\end{remark}

\section{Concluding Remarks}

We have shown that Gray's second curvature identity
holds for every invariant metric on generalized flag
manifolds with two isotropy summands, for both
invariant almost complex structures. In contrast,
the Kähler nullity depends on the metric parameters
and on the choice of almost complex structure.
Assuming
$[\mathfrak m_1,\mathfrak m_1]_{\mathfrak m_2}
\neq 0$,
the nullity for $J_2=(+,-)$ is trivial away from
$\lambda_2=2\lambda_1$ and coincides with the
integrable, totally geodesic distribution
$\mathcal D_2$ at that ratio.
For $J_+=(+,+)$, it equals the non-integrable
distribution $\mathcal D_1$ away from the
Kähler locus. Thus, Gray's second identity
alone does not imply integrability of the
Kähler nullity, highlighting the role of the
quasi-Kähler hypothesis in Gray's Theorem.

The normalized homogeneous Ricci flow reduces,
after a positive time reparametrization, to
\[
\frac{du}{d\tau}=-c(u-a)(u-2),
\qquad
a=\frac{4d_2}{d_1+2d_2},
\quad c>0.
\]
The Kähler--Einstein ratio $u=2$ is attracting,
whereas the second Einstein ratio $u=a$ is
repelling. Although Gray's second identity
is preserved along every invariant trajectory,
the Kähler nullity need not have constant
dimension under convergence. Indeed, for
$u(0)>a$ and $u(0)\neq 2$, the nullity
associated with the fixed non-integrable
structure $J_2$ remains trivial at every
finite time, while the limiting
Kähler--Einstein metric, under volume
normalization, has nullity $\mathcal D_2$.

A further direction is to investigate the
relationship between Kähler nullity and
horizontal pseudoholomorphic curves in the
homogeneous fibrations of Section 3.2.
In the twistor fibration
$\mathbb{CP}^3\to S^4$, the invariant
structure $J_2=(+,-)$ agrees with
$J_+=(+,+)$ on the horizontal distribution
$\mathcal D_1$ and differs from it by a sign
on the vertical distribution $\mathcal D_2$.
Consequently, a horizontal curve is
$J_2$-holomorphic if and only if it is
$J_+$-holomorphic. Horizontal holomorphic
curves for the standard complex structure
were studied by Bryant \cite{Bryant1982}
in connection with minimal surfaces in
$S^4$. Pseudoholomorphic curves for the
nearly Kähler structure, corresponding to
$J_2$, were investigated by Xu
\cite{Xu2010}, including horizontal and
null-torsion curves and their relationship
with holomorphic contact geometry.
Related developments concerning transverse
pseudoholomorphic curves in nearly Kähler
$\mathbb{CP}^3$ appear in
\cite{Aslan2022}.

Since the vertical distribution $\mathcal D_2$
is integrable and totally geodesic
(Theorem~3.6), and coincides with the
$J_2$-K\"ahler nullity precisely at $u=2$
(Theorem~3.7), a natural direction for
further research is to investigate how this
distinguished metric ratio interacts with
the extrinsic geometry and deformation
theory of horizontal $J_2$-holomorphic
curves. In particular, one may study the
second fundamental form of such curves
and its relationship with the horizontal
and vertical splitting
$TM=\mathcal D_1\oplus\mathcal D_2$.
The fibration
$G_2/U(2)\to G_2/SO(4)$ provides a
related setting for this investigation.

\end{document}